\documentclass[a4paper,12pt]{amsart}
\usepackage{hyperref}
\newtheorem{theorem}{Theorem}[section]
\newtheorem{lemma}[theorem]{Lemma}
\newtheorem{proposition}[theorem]{Proposition}
\newtheorem{corollary}[theorem]{Corollary}
\newtheorem{claim}{Claim}
\theoremstyle{definition}

\theoremstyle{remark}

\numberwithin{equation}{section}
\DeclareMathOperator{\Div}{div}
\DeclareMathOperator{\Span}{Span}
\DeclareMathOperator{\ls}{LS}
\DeclareMathOperator{\tr}{Tr}
\DeclareMathOperator{\norm}{Norm}
\DeclareMathOperator{\sat}{Sat}
\DeclareMathOperator{\lie}{Lie}
\DeclareMathOperator{\imm}{Im}
\newcommand{\Prj}{\mathcal{P}}
\newcommand{\Q}{\mathcal{Q}}
\newcommand{\gen}{\mathcal{L}}
\newcommand{\es}{\mathbf{E}}
\newcommand{\R}{\mathbf{R}}
\newcommand{\N}{\mathbf{N}}
\newcommand{\z}{\mathbf{Z}}
\newcommand{\q}{\mathbf{Q}}
\newcommand{\K}{\mathcal{K}}
\newcommand{\enne}{\mathcal{N}}
\newcommand{\dem}{A(\enne)}
\newcommand{\tK}{\widetilde{\mathcal{K}}}
\newcommand{\effe}{\mathcal{F}}
\newcommand{\gU}{\mathfrak{U}}
\newcommand{\gR}{\mathfrak{R}}
\newcommand{\gS}{\mathfrak{S}}
\newcommand{\e}{\mathtt{e}}
\newcommand{\im}{\mathtt{i}}
\newcommand{\ka}{\mathbf{k}}
\newcommand{\h}{\mathbf{h}}
\newcommand{\el}{\mathbf{l}}
\newcommand{\m}{\mathbf{m}}
\newcommand{\n}{\mathbf{n}}
\newcommand{\g}{\mathbf{g}}
\newcommand{\tsum}{\sideset{}{^*}\sum}
\newcommand{\set}[2]{\left\{\,#1\,|\,#2\,\right\}}
\newcommand{\dd}[1]{\delta_{#1}}
\newcommand{\diff}[1]{\frac{\partial}{\partial #1}}
\begin{document}
\title[Ergodicity of the finite dimensional...]{Ergodicity of the
finite dimensional approximation of the 3D Navier-Stokes equations
forced by a degenerate noise}
\author[M. Romito]{Marco Romito}
\address{Dipartimento di Matematica, Universit\`a di Firenze\\ Viale 
         Morgagni 67/a, 50134 Firenze, Italia}
\email{romito@math.unifi.it}
\begin{abstract}
We prove ergodicity of the finite dimensional approximations of the
three dimensional Navier-Stokes equations, driven by a random force. The
forcing noise acts only on a few modes and some algebraic conditions on
the forced modes are found that imply the ergodicity. The convergence
rate to the unique invariant measure is shown to be exponential.
\end{abstract}
\maketitle

\section{Introduction}

The uniqueness of statistical steady states for the Navier-Stokes
equations is a less known but nevertheless important problem in the
mathematical theory of turbulence. The question is completely open in
dimension three, mainly because, due to the lack of uniqueness of the
equations, there is no way yet to give meaning to the mathematical
objects involved in the subject.

In the present paper the property of ergodicity is proved for the finite
dimensional approximations of the three-dimensional Navier-Stokes
equations, driven by a random force. The same problem has
been solved in two dimensions by E and Mattingly \cite{EMa}.

Such result can have a qualitative interest for the statistical
behaviour of an incompressible fluid. Indeed, if the Kolmogorov theory
of turbulence is taken into account, one can believe that the cascade of
energy, responsible of the transport of the energy through the scales,
is effective in the inertial range so that at smaller scales only the
dissipation ends up to be relevant. Hence the long-time statistical
properties of the fluid can be sufficiently depicted by the low modes of
the velocity field. In some sense, if the ultraviolet cut-off is
sufficiently large, in order to capture all the important modes, the
corresponding invariant measure gives the real behaviour of the
fluid. In view of these considerations, the conclusions of the paper can
give also both a hint and a possible starting point for the analysis of
the infinite dimensional case.

\smallskip

We consider a finite dimensional truncation of the three dimensional
Navier-Stokes equations, driven by a random force, with periodic
boundary conditions. The proof of ergodicity is classical and it is
developed in two steps. Firstly we prove that the transition probability
densities are regular, by checking that the diffusion operator is
hypoelliptic (the H\"ormander condition). Then we show that the Markov
process is irreducible, in the sense that each open set is visited with
positive probability at each time. For this aim we study the associated
control problem (see Section \ref{control}). Irreducibility for the
infinite dimensional equations was firstly proved by Flandoli
\cite{Fla}, under the assumption that the noise acts on all modes.

Both these properties, strong Feller and irreducibility, are implied by
an algebraic condition on the set of indices corresponding to the modes
forced by the noise. The condition \textsl{essentially} means that it is
possible to obtain any index as a sum of some of the forced indices.
One can see this mechanism as a geometrical realisation of the cascade
of energy, since the non-linear term transmits the random forcing from
the few forced modes to all the other modes. As an example we show that
the algebraic condition is satisfied if the three lowest modes are
forced.

Recently, many authors have applied the techniques we have used here,
such us the hypoellipticity for degenerate diffusions, or the general
theory for Markov chains and Markov processes collected and developed by
Meyn and Tweedie (one can see their book \cite{MeTw}).  Among many
others we quote the papers by E and Mattingly \cite{EMa}, Eckmann and
Hairer \cite{EcHa}, Hairer \cite{Hai}, Rey-Bellet and Thomas
\cite{RBTh}, and some of the references therein.

\smallskip

The paper is organised as follows. In the first section the main
definitions are given, together with the statements of the main results
and an outline of their proofs. The technical computations and the
precise statement of some hypotheses are postponed in the following
sections. The aim is to give a light presentation of the main ideas,
without all the technicalities, which are then reserved to the
interested readers.

\subsection*{Acknowledgements}

The author wish to thanks R. Bianchini and S. Dolfi for the helpful
bibliographical suggestions on the control theory part in Section
\ref{control} and on the algebraic part in Section \ref{determining},
and F. Flandoli for the many helpful conversations.

\noindent This paper is dedicated to the memory of my father, who died
whilst I was writing it.

\section{The main theorem}\label{smain}

We consider the stochastic Navier-Stokes equations with additive noise
\begin{eqnarray*}
&du=(\nu\Delta u-(u\cdot\nabla)u-\nabla P)\,dt+dB_t\\
&\Div u=0,
\end{eqnarray*}
in the domain $[0,2\pi]^3$, with periodic boundary conditions, where $u$
is the velocity field and $P$ is the pressure field, and $B_t$ is a
Brownian motion. As usual, the equations are projected on the space of
divergence-free vector fields, in order to cause the pressure to
disappear from the equations. If we write the equations in the Fourier
components, we obtain the following infinite system of stochastic
differential equations
$$
du_\ka=\bigl[-\nu|\ka|^2 u_\ka
        -\im\sum_{\h+\el=\ka}(\ka\cdot u_\h)\bigl(u_\el-\frac{\ka\cdot u_\el}{|\ka|^2}\ka\bigr)\bigr]\,dt
        +q_\ka\,d\beta^\ka_t,\qquad \ka\in\z^3
$$
with the constraint $u_\ka\cdot\ka=0$ (it comes from the divergence-free
condition). We have made some simplifying assumptions on the noise: we assume
that the noise takes values in the space of divergence-free vector fields
and that the covariance is diagonal in the Fourier components (the
assumptions will be stated more clearly in Section \ref{noiseass}). 

In order to state the problem of the finite dimensional approximation,
fix a threshold $N$ and consider the finite subset of indices
$\K_N=\set{\ka\in\z^3}{|\ka|\le N,\,\ka\neq(0,0,0)}$. The finite
dimensional system obtained is the following
\begin{equation}\label{finiteq}
du_\ka=\bigl[-\nu|\ka|^2 u_\ka
        -\im\sum_{\substack{\h,\el\in\K_N\\\h+\el=\ka}}(\ka\cdot
u_\h)\bigl(u_\el-\frac{\ka\cdot u_\el}{|\ka|^2}\ka\bigr)\bigr]\,dt
        +q_\ka\,d\beta^\ka_t,
\qquad \ka\in\K_N
\end{equation}
(a formal derivation is given in Section \ref{fourier}). We will use the
real variables $r_\ka$, $s_\ka\in\R^3$, where $u_\ka=r_\ka+\im s_\ka$,
rather than the complex variables $u_\ka$, so that the equations are
briefly written as
$$
\left\{\begin{array}{l}
dr_\ka^i=F_{r^i_\ka}(r,s)\,dt+q^r_\ka\,d\beta^\ka_t,\\
ds_\ka^i=F_{s^i_\ka}(r,s)\,dt+q^s_\ka\,d\beta^\ka_t,
\end{array}\right.
\qquad \ka\in\K_N,\quad i=1,2,3
$$
and $q_\ka=q_\ka^r+\im q_\ka^s$. Since $u_{-\ka}=\overline{u_\ka}$, the
set of indices $\K_N$ is redundant, hence we take a smaller set $\tK$,
which takes into account the symmetries.

The solution $(r(t),s(t))$ of the above stochastic equations is a Markov
process on the state space
$$
U=\bigoplus_{\ka\in\tK}(R_\ka\oplus S_\ka),
$$
where $R_\ka$ and $S_\ka$ enclose the divergence-free condition
$r_\ka\cdot\ka=s_\ka\cdot\ka=0$ (see also \eqref{globspace} and the
following formulas). We denote by $P_t$ the transition semigroup
$$
P_t\varphi(r_0,s_0)=\es_{(r_0,s_0)}[\varphi(r(t),s(t))]
$$
with generator
\begin{equation}\label{markovgen}
\gen=F_0+\frac12\sum_{\substack{\ka\in\tK\\ i=1,2,3}}({X^r_{\ka,i}}^2+{X^s_{\ka,i}}^2)
\end{equation}
where
\begin{equation}\label{campi}
F_0=\sum_{\ka\in\tK}\sum_{i=1}^3F_{r^i_\ka}\diff{r^i_\ka}+F_{s^i_\ka}\diff{s^i_\ka},
\qquad X^r_{\ka,i}=\sum_{j=1}^3 q^r_{\ka,ij}\diff{r^j_\ka},
\qquad X^s_{\ka,i}=\sum_{j=1}^3 q^s_{\ka,ij}\diff{s^j_\ka},
\end{equation}
and by $P_t((r,s),\cdot)$ the transition probability.

The main assumption we take on the noise is that it acts on a small set
of modes. We consider the set $\enne$ of indices whose corresponding
Fourier components are forced by the noise. We assume that $\enne$ is a
\textsl{determining set of indices}, as defined in Section
\ref{determining}, which \textsl{essentially} means that each index in
$\K_N$ can be obtained as the sum of elements of $\enne$. In other words
$\enne$ should be an algebraic system of generators of $\z^3$. In
Section \ref{determining} we will give some heuristic justifications
to such claim. As a \textsl{working example}, Proposition \ref{workex}
shows that any set $\enne$ containing the three indices $(1,0,0)$,
$(0,1,0)$ and $(0,0,1)$ is a determining set of indices.

Here we are interested in stating the main result of the paper, namely the
ergodicity of the finite dimensional approximation \eqref{finiteq}

\begin{theorem}\label{main}
Assume that the Brownian motion $B_t$ satisfies the assumptions in
Section \ref{noiseass} and that the set $\enne$ defined above is a
determining set of indices. Then the system \eqref{finiteq} admits a
unique invariant measure.

Moreover, the unique invariant measure is supported on the whole state
space or, in other words, it gives positive mass to each open set.
\end{theorem}

\begin{proof}
First, we prove the existence of the invariant measure. The method is
classical and based on the Krylov-Bogoliubov method (see for example
Theorem 3.1.1 of Da Prato and Zabczyk \cite{DaPZa}). The compactness
follows by the following argument. Let
$\|u\|^2=\sum_{\ka\in\tK}|u_\ka|^2$, then by It\^o formula (using also
the first property of Lemma \ref{liapunov}),
\begin{eqnarray*}
d\|u(t)\|^2
&=& \sum_{\ka\in\tK}\bigl(2\overline{u_\ka}\cdot
F_{\ka}(u)+\tr(q_\ka^T\cdot\overline{q_\ka})\bigr)\,dt
   +2\sum_{\ka\in\tK}\overline{u_\ka}\cdot q_\ka\,d\beta^\ka_t\\
&=&-2\nu\sum_{\ka\in\tK}|\ka|^2|u_\ka|^2\,dt
   +2\sum_{\ka\in\tK}\overline{u_\ka}\cdot q_\ka\,d\beta^\ka_t
   +\sigma^2\,dt,
\end{eqnarray*}
where $\sigma^2$ is the variance of the Brownian motion $B_t$, so that
$$
\es\|u(t)\|^2+2\nu\int^t_0\|u(s)\|^2\,ds\le\es\|u(0)\|+\sigma^2t
$$
and by Gronwall lemma $\es\|u(t)\|\le \es\|u(0)\|+\frac{\sigma^2}{2\nu}$.

Uniqueness of the invariant measure is proved by means of the Doob
uniqueness theorem (see for example Theorem 4.2.1 of Da Prato and
Zabczyk \cite{DaPZa}). We just need to show that the transition
semigroup generated by the dynamics \eqref{finiteq} is strongly Feller
and irreducible.

A Markov semigroup $P_t$ is strongly Feller if $P_t\varphi$ is bounded
continuous in time and space when $\varphi$ is bounded measurable. By a
theorem of Stroock \cite{Str}, the transition semigroup is strongly
Feller if the \textsl{H\"ormander condition} holds: the Lie algebra
generated by the vector fields in \eqref{campi}, evaluated at each
point, is the state space $U$. Since $\enne$ is a determining set of
indices, from Lemma \ref{combi} it follows that the constant vector
fields of the generated Lie algebra span $U$.

A Markov semigroup is irreducible if it gives positive mass to any open
set for each initial condition and each time. It is well known (see
Stroock and Varadhan \cite{StVa}) that irreducibility is true if the
control problem (see equations \eqref{controleq}) associated to problem
\eqref{finiteq} is controllable. The last statement follows from Theorem
\ref{controllo}.

Finally, the irreducibility property implies also that the support of
the invariant measure is the whole state space.
\end{proof}

The next theorem shows, by means of general techniques developed in Meyn
and Tweedie \cite{MeTw2}, \cite{MeTw3}, that the finite approximation of
Navier-Stokes equations has good dissipation properties, strong enough
to ensure the exponential mixing of the dynamics given by the Markov
process. In order to state the result, define, for any measurable
function $f\ge 1$ and any signed measure $\mu$ on the Borel sets of $U$,
$$
\|\mu\|_f=\sup_{|g|\le f}\Big|\int g(x)\mu(dx)\Big|,
$$
and set
$$
V(r,s)=\sum_{\ka\in\tK}\sum_{i=1,2,3}({r^i_\ka}^2+{s^i_\ka}^2),\qquad
(r,s)\in U.
$$

\begin{theorem}\label{exporate}
Under the assumptions of the previous theorem, let $\pi$ be the unique
invariant measure. Then there are positive constants $C$ and $\rho$ such
that for each initial condition $(r_0,s_0)\in U$,
$$
\|P_t((r_0,s_0),\cdot)-\pi\|_f\le C\e^{-\rho t}\bigl(1+V(r_0,s_0)+\frac{\sigma^2}{2\nu}\bigr),\qquad t>0,
$$
where $f=1+V$.
\end{theorem}

\section{The Navier-Stokes equations in the Fourier
coordinates}\label{fourier}

In this section we derive the equations of the finite dimensional
approximations of the stochastic Navier-Stokes equations, with additive
noise,
\begin{eqnarray*}
&du=(\nu\Delta u-(u\cdot\nabla)u-\nabla P)\,dt+dB_t\\
&\Div u=0,
\end{eqnarray*}
in the domain $[0,2\pi]^3$, with periodic boundary conditions,
in the Fourier components.

Consider the Fourier basis $(\e^{\im\ka\cdot x})_{\ka\in\z^3}$ of
$L^2([0,2\pi]^3)$. First, assume that the applied random force has zero
average, so that the centre of mass of the fluid moves with constant
velocity. Hence, without loss of generality, we can assume that
$$
u_{\bf 0}=P_{\bf 0}=0.
$$
The projection onto the space of divergence-free vector fields is
defined as
$$
\Prj(a\e^{\im\ka\cdot x})
=\big(a-\frac{\ka\otimes\ka}{|\ka|^2}\cdot a\big)\e^{\im\ka\cdot x}
=\big(a-\frac{a\cdot\ka}{|\ka|^2}\ka\big)\e^{\im\ka\cdot x},
$$
where $|\cdot|$ is the Euclidean norm in $\R^3$. Notice that
$$
\Div u=0\quad\text{means}\quad \ka\cdot u_\ka=0\quad \text{for each
}\ka.
$$

\subsection{Assumptions on the noise}\label{noiseass}

For the sake of simplicity, some simplifying assumptions will be
done. First we assume that the covariance $\Q$ of the noise is diagonal
in the Fourier basis, so that we can write
$$
\Q v=\sum_{\ka\in\z^3}(q_\ka\cdot v_\ka)\e^{\im\ka\cdot x}.
$$
Moreover we assume that $q_\ka^T\cdot\ka=\mathbf{0}$ for each $\ka$,
this implies that the Brownian motion takes values in the space of
divergence-free vector fields. The Brownian motion has finite variance
that we denote by $\sigma^2$. We assume also that for each $\ka$, the
real and the imaginary parts of the $3\times3$ matrix $q_\ka$, if not
zero, have rank $2$. This is an assumption \textsl{in the small} of
non-degeneracy, since we ask that, if a mode is forced, it is fully
forced in its $4$ components. As a first consequence of our assumptions,
the operators $\Q$ and $\Prj$ commute.

The main assumption of the paper is that the noise acts only on a few
components, namely most of the matrices $q_\ka$ are zero. We define the
set $\enne\subset\z^3$ of stochastically forced indices, that is the set
of $\ka$s such that $q_\ka\not\equiv0$.

\subsection{The equation in the Fourier modes}
We write
$$
u(t,x)=\sum_{\ka\in\z^3}u_\ka(t)\e^{\im\ka\cdot x}
$$
and, by means of the operator $\Prj$, we project the equations in the
space of divergence-free vector fields, so that the pressure
disappears. We obtain the following infinite system of stochastic
differential equations (see also Gallavotti \cite{Gal}, Chapter 2, where
the author gives also a interpretation of the physics of the fluid in
terms of the Fourier coordinates)
\begin{eqnarray*}
&&du_\ka=\bigl[-\nu|\ka|^2 u_\ka
          -\im\sum_{\substack{\h,\el\in\z^3\\\h+\el=\ka}}(\ka\cdot
u_\h)\bigl(u_\el-\frac{\ka\cdot u_\el}{|\ka|^2}\ka\bigr)\bigr]\,dt
          +q_\ka\,d\beta^\ka_t,\\
&&u_\ka\cdot\ka=0
\end{eqnarray*}
where $(\beta^\ka_t)_{t\ge0}$ are independent three-dimensional Brownian
motions, and the nonlinear term has been obtained in the following way:
\begin{eqnarray*}
\Prj(u\cdot\nabla)u
&=& \im\Prj\sum_{\ka\in\z^3}\sum_{\h+\el=\ka}(\el\cdot
u_\h)u_\el\e^{\im\ka\cdot x}\\
&=& \im\sum_{\ka\in\z^3}\sum_{\h+\el=\ka}(\el\cdot
u_\h)\big(u_\el-\frac{\ka\cdot u_\el}{|\ka|^2}\ka\big)\e^{\im\ka\cdot x}\\
&=& \im\sum_{\ka\in\z^3}\sum_{\h+\el=\ka}(\ka\cdot
u_\h)\big(u_\el-\frac{\ka\cdot u_\el}{|\ka|^2}\ka\big)\e^{\im\ka\cdot x}.
\end{eqnarray*}

\subsection{The finite dimensional approximation}\label{finricava}

Let $N\in\N$ and set
$$
\K_N=\set{\ka\in\z^3}{\ka\neq(0,0,0),\ |\ka|_\infty\le N}.
$$
where $|\cdot|_\infty$ is the sup-norm in $\R^3$.
We project the equation in the space spanned by $(\e^{\im\ka\cdot
x})_{\ka\in\K_N}$, with coefficients in $\R^3$,  and for this aim we set
$$
u(t,x)=\sum_{\ka\in\K_N}u_\ka\e^{\im\ka\cdot x}.
$$
The equation in the finite dimensional approximation is
$$
du_\ka=\bigl[-\nu|\ka|^2 u_\ka
        -\im\sum_{\substack{\h,\el\in\K_N\\\h+\el=\ka}}(\ka\cdot
u_\h)\bigl(u_\el-\frac{\ka\cdot u_\el}{|\ka|^2}\ka\bigr)\bigr]\,dt
        +q_\ka\,d\beta^\ka_t,
\qquad \ka\in\K_N.
$$

We set
$$
u_\ka=(r_\ka^j+\im s_\ka^j)_{j=1,2,3},
$$
where $\ka\cdot r_\ka=\ka\cdot s_\ka=0$ and $r_\ka^j$, $s_\ka^j$,
$j=1,2,3$, are real-valued. Since $u_{-\ka}=\overline{u_\ka}$, we are
going to choose a smaller set of indices $\ka\in\K_N$ in order to
take into account that some equations in the system are redundant. We set 
\begin{eqnarray*}
\K^1_N&=&\set{\ka\in\z^3}{|\ka|_\infty\le N,\ k_3>0}\\
\K^2_N&=&\set{\ka\in\z^3}{|\ka|_\infty\le N,\ k_3=0,\ k_2>0}\\
\K^3_N&=&\set{\ka\in\z^3}{|\ka|_\infty\le N,\ k_3=k_2=0,\ k_1>0}\\
\end{eqnarray*}
and
$$
\tK=\K_N^1\cup\K_N^2\cup\K_N^3,
$$
in such a way that
$$
\K_N=\tK\cup(-\tK)\quad\text{and}\quad \tK\cap(-\tK)=\emptyset.
$$
Notice that $\#(\tK)=\frac12[(2N+1)^3-1]$, we call such number $D$.
Now, if $\ka\in\tK$, the sum extended to all pairs of indices $\h$,
$\el$ such that $\h+\el=\ka$ can be written in the following way:
$$
\sum_{\substack{\h+\el=\ka\\ \h,\el\in\K_N}}
=\sum_{\substack{\h+\el=\ka\\ \h,\el\in\tK}}
+\sum_{\substack{\h+\el=\ka\\ \h\in\tK\\ \el\in-\tK}}
+\sum_{\substack{\h+\el=\ka\\ \h\in-\tK\\ \el\in\tK}},
$$
since if $\h$, $\el\not\in\tK$, $\ka$ does not belong to $\tK$ as
well. We denote by $\tsum$ the sum extended to indices in $\tK$. With
this position
\begin{eqnarray*}
\sum_{\substack{\h+\el=\ka\\\h,\el\in\K_N}}(\ka\cdot u_\h)\big(u_\el-\frac{\ka\cdot
u_\el}{|\ka|^2}\ka\big)
&=& \tsum_{\h+\el=\ka}(\ka\cdot u_\h)\big(u_\el-\frac{\ka\cdot
u_\el}{|\ka|^2}\ka\big)
   +\tsum_{\h-\el=\ka}(\ka\cdot u_\h)\big(\overline{u_\el}-\frac{\ka\cdot
\overline{u_\el}}{|\ka|^2}\ka\big)\\
& & +\tsum_{\el-\h=\ka}(\ka\cdot\overline{u_\h})\big(u_\el-\frac{\ka\cdot
u_\el}{|\ka|^2}\ka\big)\\
\end{eqnarray*}
so that the equations become
\begin{eqnarray*}
du_\ka
+\Bigl(\nu|\ka|^2u_\ka
 + \im\tsum_{\h+\el=\ka}(\ka\cdot u_\h)\big(u_\el-\frac{\ka\cdot
   u_\el}{|\ka|^2}\ka\big)
&+&\im\tsum_{\h-\el=\ka}(\ka\cdot u_\h)\big(\overline{u_\el}-\frac{\ka\cdot
   \overline{u_\el}}{|\ka|^2}\ka\big)+\\
&+&\im\tsum_{\el-\h=\ka}(\ka\cdot\overline{u_\h})\big(u_\el-\frac{\ka\cdot
   u_\el}{|\ka|^2}\ka\big)\Bigr)\,dt
=q_\ka\,d\beta^\ka_t.
\end{eqnarray*}
It is convenient to write explicitly the equations relative to the real
and imaginary part of $u_\ka$,
\begin{eqnarray*}
dr_\ka
+\Bigl(\nu|\ka|^2r_\ka
&-&\tsum_{\h+\el=\ka}(\ka\cdot r_\h)\big(s_\el-\frac{\ka\cdot
   s_\el}{|\ka|^2}\ka\big)+(\ka\cdot s_\h)\big(r_\el-\frac{\ka\cdot
   r_\el}{|\ka|^2}\ka\big)\\
&+&\tsum_{\h-\el=\ka}(\ka\cdot r_\h)\big(s_\el-\frac{\ka\cdot
   s_\el}{|\ka|^2}\ka\big)-(\ka\cdot s_\h)\big(r_\el-\frac{\ka\cdot
   r_\el}{|\ka|^2}\ka\big)\\
&-&\tsum_{\el-\h=\ka}(\ka\cdot r_\h)\big(s_\el-\frac{\ka\cdot
   s_\el}{|\ka|^2}\ka\big)-(\ka\cdot s_\h)\big(r_\el-\frac{\ka\cdot
   r_\el}{|\ka|^2}\ka\big)\Bigr)\,dt=q_\ka^rd\beta^\ka_t\\
\end{eqnarray*}
and
\begin{eqnarray*}
ds_\ka
+\Bigl(\nu|\ka|^2s_\ka
&+&\tsum_{\h+\el=\ka}(\ka\cdot r_\h)\big(r_\el-\frac{\ka\cdot
   r_\el}{|\ka|^2}\ka\big)-(\ka\cdot s_\h)\big(s_\el-\frac{\ka\cdot
   s_\el}{|\ka|^2}\ka\big)\\
&+&\tsum_{\h-\el=\ka}(\ka\cdot r_\h)\big(r_\el-\frac{\ka\cdot
   r_\el}{|\ka|^2}\ka\big)+(\ka\cdot s_\h)\big(s_\el-\frac{\ka\cdot
   s_\el}{|\ka|^2}\ka\big)\\
&+&\tsum_{\el-\h=\ka}(\ka\cdot r_\h)\big(r_\el-\frac{\ka\cdot
   r_\el}{|\ka|^2}\ka\big)+(\ka\cdot s_\h)\big(s_\el-\frac{\ka\cdot
   s_\el}{|\ka|^2}\ka\big)\Bigr)\,dt=q_\ka^sd\beta^\ka_t.\\
\end{eqnarray*}
In view of the above formulas, we set
\begin{eqnarray}\label{fr}
F_{r^i_\ka}=
-\nu|\ka|^2r^i_\ka
&+&\tsum_{\h+\el=\ka}(\ka\cdot r_\h)\big(s^i_\el-\frac{\ka\cdot
   s_\el}{|\ka|^2}k_i\big)+(\ka\cdot s_\h)\big(r^i_\el-\frac{\ka\cdot
   r_\el}{|\ka|^2}k_i\big)\notag\\
&-&\tsum_{\h-\el=\ka}(\ka\cdot r_\h)\big(s^i_\el-\frac{\ka\cdot
   s_\el}{|\ka|^2}k_i\big)-(\ka\cdot s_\h)\big(r^i_\el-\frac{\ka\cdot
   r_\el}{|\ka|^2}k_i\big)\\
&+&\tsum_{\el-\h=\ka}(\ka\cdot r_\h)\big(s^i_\el-\frac{\ka\cdot
   s_\el}{|\ka|^2}k_i\big)-(\ka\cdot s_\h)\big(r^i_\el-\frac{\ka\cdot
   r_\el}{|\ka|^2}k_i\big)\notag
\end{eqnarray}
and
\begin{eqnarray}\label{fs}
F_{s^i_\ka}=
-\nu|\ka|^2s^i_\ka
&-&\tsum_{\h+\el=\ka}(\ka\cdot r_\h)\big(r^i_\el-\frac{\ka\cdot
   r_\el}{|\ka|^2}k_i\big)-(\ka\cdot s_\h)\big(s^i_\el-\frac{\ka\cdot
   s_\el}{|\ka|^2}k_i\big)\notag\\
&-&\tsum_{\h-\el=\ka}(\ka\cdot r_\h)\big(r^i_\el-\frac{\ka\cdot
   r_\el}{|\ka|^2}k_i\big)+(\ka\cdot s_\h)\big(s^i_\el-\frac{\ka\cdot
   s_\el}{|\ka|^2}k_i\big)\\
&-&\tsum_{\el-\h=\ka}(\ka\cdot r_\h)\big(r^i_\el-\frac{\ka\cdot
   r_\el}{|\ka|^2}k_i\big)+(\ka\cdot s_\h)\big(s^i_\el-\frac{\ka\cdot
   s_\el}{|\ka|^2}k_i\big).\notag
\end{eqnarray}

\section{The Lie Algebra generated by the dynamics}

The state space of the Markov process $(r(t),s(t))$ which is solution of
the equations stated above is a linear space $U\subset\R^{6D}$, where
$D=\#\tK$, given by
\begin{equation}\label{globspace}
U=\bigoplus_{\ka\in\tK}(R_\ka\oplus S_\ka),
\end{equation}
and each element of $U$ is labelled $(r,s)$, with
$r=(r_\ka^1,r_\ka^2,r_\ka^3)_{\ka\in\tK}$ and
$s=(s_\ka^1,s_\ka^2,s_\ka^3)_{\ka\in\tK}$, and
\begin{eqnarray*}
R_\ka&=&\set{(r,s)\in\R^{6D}}{r_\ka\cdot\ka=0,\ s_\ka=0,\ r_\h=s_\h=0,\ \h\neq\ka}\\
S_\ka&=&\set{(r,s)\in\R^{6D}}{s_\ka\cdot\ka=0,\ r_\ka=0,\ r_\h=s_\h=0,\ \h\neq\ka}.
\end{eqnarray*}

In the same way, we can define the Lie algebra $\gU$ corresponding to the
vector space $U$,
\begin{equation}\label{globLie}
\gU=\Bigl\{G\,\Big|\,G=\sum_{\substack{\ka\in\tK\\ i=1,2,3}}
            G_{r^i_\ka}\frac\partial{\partial r^i_\ka}
           +G_{s^i_\ka}\frac\partial{\partial s^i_\ka}
    \text{ and }\ka\cdot G_{r_\ka}=\ka\cdot G_{s_\ka}=0\,\Bigr\}.
\end{equation}
We define also the subspaces $\gU_\ka=\gR_\ka\oplus\gS_\ka$ of $\gU$ of
constant vector fields, where
$$
\gR_\ka=\big\{\sum_{i=1,2,3}r^i_\ka\diff{r^i_\ka}\,|\,r_\ka\in R_\ka\,\big\}
\quad\text{and}\quad
\gS_\ka=\big\{\sum_{i=1,2,3}s^i_\ka\diff{s^i_\ka}\,|\,s_\ka\in S_\ka\,\big\}
$$

In this section, we want to find some reasonable conditions on the set
$\enne$ of forced modes (such set has been defined in Section
\ref{noiseass}) in such a way that the algebra generated by the fields
\begin{equation}\label{generatori}
\{F_0\}\cup\gU_\ka\qquad\ka\in\enne,
\end{equation}
where
$$
F_0=\sum_{\substack{\ka\in\tK\\ i=1,2,3}}
 F_{r^i_\ka}\diff{r^i_\ka}
+F_{s^i_\ka}\diff{s^i_\ka},
$$
and $F_{r^i_\ka}$ and $F_{s^i_\ka}$ have been defined respectively in
\eqref{fr} and \eqref{fs}, contains all the constant vector fields of
$\gU$. In particular, it follows that the H\"ormander condition holds,
that is the generated Lie algebra, evaluated at each point of $U$,
gives $U$ itself. We start with some computations that will be useful in
the sequel.

\begin{lemma}\label{contacci}
Let $\m$, $\n\in\tK$ and $V\in\gU_\m$, $W\in\gU_\n$, with
$$
V=\sum_{j=1}^3v^r_j\diff{r^j_\m}+v^s_j\diff{s^j_\m},
\qquad
W=\sum_{j=1}^3w^r_j\diff{r^j_\m}+w^s_j\diff{s^j_\m},
$$
then
\begin{enumerate}
\item[\textit{(i)}] if $\ka=\m+\n$, $\h=\n-\m$ and $\g=\m-\n$,
\begin{eqnarray*}
[[F_0,V],W]
&=&\bigl[(v^r\cdot\ka)P_\ka(w^s)+(w^s\cdot\ka)P_\ka(v^r)+(v^s\cdot\ka)P_\ka(w^r)+(w^r\cdot\ka)P_\ka(v^s)\bigr]\cdot\diff{r_\ka}\\
& &+\bigl[(v^s\cdot\ka)P_\ka(w^s)+(w^s\cdot\ka)P_\ka(v^s)-(v^r\cdot\ka)P_\ka(w^r)-(w^r\cdot\ka)P_\ka(v^r)\bigr]\cdot\diff{s_\ka}\\
& &+\bigl[(v^r\cdot\h)P_\h(w^s)+(w^s\cdot\h)P_\h(v^r)-(v^s\cdot\h)P_\h(w^r)-(w^r\cdot\h)P_\h(v^s)\bigr]\cdot\diff{r_\h}\\
& &-\bigl[(v^r\cdot\h)P_\h(w^r)+(w^r\cdot\h)P_\h(v^r)+(v^s\cdot\h)P_\h(w^s)+(w^s\cdot\h)P_\h(v^s)\bigr]\cdot\diff{s_\h}\\
& &+\bigl[(v^s\cdot\g)P_\g(w^r)+(w^r\cdot\g)P_\g(v^s)-(v^r\cdot\g)P_\g(w^s)-(w^s\cdot\g)P_\g(v^r)\bigr]\cdot\diff{r_\g}\\
& &-\bigl[(v^r\cdot\g)P_\g(w^r)+(w^r\cdot\g)P_\g(v^r)+(v^s\cdot\g)P_\g(w^s)+(w^s\cdot\g)P_\g(v^s)\bigr]\cdot\diff{s_\g},
\end{eqnarray*}
where $P_\ka$ is the projection of $\R^3$ on the plane orthogonal to the
vector $\ka$, and in the above formula the terms corresponding to
indices out of $\tK$ are zero;
\item[\textit{(ii)}] if there is $q\in\q$ such that $\n=q\,\m$, then $[[F_0,V],W]=0,$
\item[\textit{(iii)}] $[[F_0,V],W]=\frac12[[F_0,V+W],V+W]$.
\end{enumerate}
\end{lemma}

\begin{proof}
We compute the derivatives of the components of $F_0$ (defined in
\eqref{fr} and \eqref{fs}),
\begin{eqnarray*}
\frac{\partial F_{r^i_\ka}}{\partial r^j_\m}&=&
-\nu|\ka|^2\dd{ij}\dd{\ka\m}
+k_j(s^i_{\ka-\m}-s^i_{\m-\ka}+s^i_{\m+\ka})
+\ka\cdot(s_{\ka-\m}-s_{\m-\ka}+s_{\m+\ka})\big(\dd{ij}-\frac{2k_ik_j}{|\ka|^2}\big)\\
\frac{\partial F_{r^i_\ka}}{\partial s^j_\m}&=&
k_j(r^i_{\ka-\m}+r^i_{\m-\ka}-r^i_{\m+\ka})
+\ka\cdot(r_{\ka-\m}+r_{\m-\ka}-r_{\m+\ka})\big(\dd{ij}-\frac{2k_ik_j}{|\ka|^2}\big)\\
\frac{\partial F_{s^i_\ka}}{\partial r^j_\m}&=&
-k_j(r^i_{\ka-\m}+r^i_{\m-\ka}+r^i_{\m+\ka})
-\ka\cdot(r_{\ka-\m}+r_{\m-\ka}+r_{\m+\ka})\big(\dd{ij}-\frac{2k_ik_j}{|\ka|^2}\big)\\
\frac{\partial F_{s^i_\ka}}{\partial s^j_\m}&=&
-\nu|\ka|^2\dd{ij}\dd{\ka\m}
+k_j(s^i_{\ka-\m}-s^i_{\m-\ka}-s^i_{\m+\ka})
+\ka\cdot(s_{\ka-\m}-s_{\m-\ka}-s_{\m+\ka})\big(\dd{ij}-\frac{2k_ik_j}{|\ka|^2}\big)
\end{eqnarray*}
and the second derivatives (we have set
$A^i_{jl}(\ka)=\dd{il}k_j+\dd{ij}k_l-2\frac{k_ik_jk_l}{|\ka|^2}$),
$$
\frac{\partial^2 F_{r^i_\ka}}{\partial r^l_\n\partial r^j_\m}=
\frac{\partial^2 F_{r^i_\ka}}{\partial s^l_\n\partial s^j_\m}=0,
\qquad
\frac{\partial^2 F_{r^i_\ka}}{\partial s^l_\n\partial r^j_\m}=
(\dd{\n,\ka-\m}-\dd{\n,\m-\ka}+\dd{\n,\m+\ka})A^i_{ij}(\ka)
$$
and
\begin{eqnarray*}
\frac{\partial^2 F_{s^i_\ka}}{\partial r^l_\n\partial r^j_\m}&=&
-(\dd{\n,\ka-\m}+\dd{\n,\m-\ka}+\dd{\n,\m+\ka})A^i_{jl}(\ka),\qquad
\frac{\partial^2 F_{s^i_\ka}}{\partial s^l_\n\partial r^j_\m}=0\\
\frac{\partial^2 F_{s^i_\ka}}{\partial s^l_\n\partial s^j_\m}&=&
(\dd{\n,\ka-\m}-\dd{\n,\m-\ka}-\dd{\n,\m+\ka})A^i_{jl}(\ka)
\end{eqnarray*}
with the agreement that everything concerning indices out of the set
$\tK$ is zero.

Let $V\in\gU_\m$ and $W\in\gU_\n$ as in the statement of the lemma, then
by computing the bracket we obtain
$$
[[F_0,V],W]=\sum_{\ka\in\tK}\sum_{i,j,l=1}^3
\bigl(v^s_jw_l^r\frac{\partial^2F_{r^i_\ka}}{\partial s^j_\m\partial
r^l_\n}+v^r_jw_l^s\frac{\partial^2F_{r^i_\ka}}{\partial r^j_\m\partial
s^l_\n}\bigr)\diff{r^i_\ka}
+\bigl(v^r_jw_l^r\frac{\partial^2F_{s^i_\ka}}{\partial r^j_\m\partial
r^l_\n}+v^s_jw_l^s\frac{\partial^2F_{s^i_\ka}}{\partial s^j_\m\partial
s^l_\n}\bigr)\diff{s^i_\ka}.
$$
We analyse the coefficients of the $\partial_{r^i_\ka}$-components:
\begin{eqnarray*}
\lefteqn{v^s_jw_l^r\frac{\partial^2F_{r^i_\ka}}{\partial s^j_\m\partial
r^l_\n}+v^r_jw_l^s\frac{\partial^2F_{r^i_\ka}}{\partial r^j_\m\partial
s^l_\n}=}\\
&=&\sum_{j,l=1}^3(\dd{\m,\ka-\n}-\dd{\m,\n-\ka}+\dd{\m,\n+\ka})A^i_{jl}(\ka)v_j^sw_l^r
                +(\dd{\n,\ka-\m}-\dd{\n,\m-\ka}+\dd{\n,\m+\ka})A^i_{jl}(\ka)v_j^rw_l^s\\
&=&(\dd{\m,\ka-\n}-\dd{\m,\n-\ka}+\dd{\m,\n+\ka})
      [(v^s\cdot\ka)P_\ka(w^r)_i+(w^r\cdot\ka)P_\ka(v^s)_i]\\
& &+(\dd{\n,\ka-\m}-\dd{\n,\m-\ka}+\dd{\n,\m+\ka})
      [(v^r\cdot\ka)P_\ka(w^s)_i+(w^s\cdot\ka)P_\ka(v^r)_i],
\end{eqnarray*}
where $P_\ka(v)_i=v_i-\frac{k_i}{|\ka|^2}(v\cdot\ka)$. In a similar way
it is possible to treat the coefficients of the
$\partial_{s^i_\ka}$-components, and claim \textit{(i)} is true.

If $\n=q\,\m$, it follows that
$$
v^r\cdot\ka=v^s\cdot\ka=w^r\cdot\ka=w^s\cdot\ka=0
$$
with $\ka=\m+\n$, $\m-\n$ and $\n-\m$, so that using property
\textit{(i)} of this lemma, claim \textit{(ii)} holds true. Finally, if
$V\in\gU_\m$ and $W\in\gU_\n$, by the Jacobi identity,
$$
[[F_0,V+W],V+W]
=[[F_0,V],V]+[[F_0,V],W]+[[F_0,W],V]+[[F_0,W],W]
=2[[F_0,V],W].
$$
\end{proof}

The computations of the above lemma show that the non-linear term mixes
and combines the components. In some sense, this mechanism can be
considered as a geometrical counterpart of the cascade of energy. Our
aim is to understand for which sets $\enne$ of forced modes the
evaluation of the Lie algebra generated by the fields
\eqref{generatori}, gives $U$.  We define the set $\dem\subset\K_N$ of
indices $\ka\in K_N$ such that the constant vector fields corresponding
to $\ka$ (or to $-\ka$, depending on $\ka\in\tK$ or $-\ka\in\tK$) are in
the Lie algebra generated by the vector fields
\eqref{generatori}. Obviously, $\enne\subset\dem$, and our aim is to
show that $\dem=\K_N$.

\begin{lemma}\label{combi}
Let $\enne$ be a subset of indices and define the set $\dem$ as above.
\begin{enumerate}
\item[\textit{(i)}] If $\m\in\dem$, then also $-\m\in\dem$,
\item[\textit{(ii)}] if $\m$, $\n$ are in $\dem$, $\m+\n$ is in
$\K_N$, $\m$ and $\n$ are linearly independent and $|\m|\neq|\n|$, then
$\m+\n\in\dem$.
\end{enumerate}
\end{lemma}
\begin{proof}
The first property follows from the fact that $u_{-\m}=\overline{u_{\m}}$.
In order to show the second claim, take $\m$ and $\n$ in $\dem\cap\tK$
and assume that $\ka=\m+\n\in\tK$. The claim follows if $\m+\n\in\dem$.

Let
$$
V^r=\sum_{i=1}^3v_i\diff{r^i_\m},\qquad
V^s=\sum_{i=1}^3v_i\diff{s^i_\m},\qquad
W^r=\sum_{i=1}^3w_i\diff{r^i_\n},\qquad
W^s=\sum_{i=1}^3w_i\diff{s^i_\n},
$$
with $v\cdot\m=w\cdot\n=0$. Then, by property \textit{(i)} of the
previous lemma,
\begin{eqnarray*}
[[F_0,V^r],W^s]+[[F_0,V^s],W^r]
&=&2\bigl((v\cdot\ka)P_\ka(w)+(w\cdot\ka)P_\ka(v)\bigr)\cdot\diff{r_\ka},\\
\,[[F_0,V^r],W^r]-[[F_0,V^s],W^s]
&=&2\bigl((v\cdot\ka)P_\ka(w)+(w\cdot\ka)P_\ka(v)\bigr)\cdot\diff{s_\ka}.
\end{eqnarray*}
Now, let $M$, $E$ two vectors in $\R^3$ such that $\{\ka,M,E\}$ is a
basis of $\R^3$, $M$ and $E$ span $\set{x\in\R^3}{x\cdot\ka=0}$ and
$\m$, $\n$ are in $\Span[\ka,M]$. Choose
$$
v=\lambda_1\ka+\mu_1 M+\nu_1 E,\qquad
w=\lambda_2\ka+\mu_2 M+\nu_2 E,
$$
then, by the assumptions on $\m$ and $\n$, it is always possible to
choose the coefficients $\lambda_1$, $\mu_1$, $\nu_1$, $\lambda_2$,
$\mu_2$, $\nu_2$ in such a way that
$(v\cdot\ka)P_\ka(w)+(w\cdot\ka)P_\ka(v)$ is any vector in
$\Span[M,E]$. In other words, $\gU_\ka$ is contained in the Lie algebra
generated by the vector fields \eqref{generatori}. In the same way, if
$\h=\n-\m\in\tK$ (or if $\g=\m-\n\in\tK$), the conclusion follows by
taking $[[F_0,V^r],W^s]-[[F_0,V^s],W^r]$ and
$[[F_0,V^r],W^r]+[[F_0,V^s],W^s]$.
\end{proof}

\section{Determining sets of indices}\label{determining}

In view of Lemma \ref{combi}, we call a subset $\enne$ of $\K_N$ a
\textsl{determining set of indices} for the ultraviolet cut-off $N$, if
$\enne$ generates the cube $K_N$ in the sense that $\dem=\K_N$, where
$\dem$ has been defined in the previous section. Lemma \ref{combi} shows
us which is the algebraic structure of such set. Namely, $\dem$ is
symmetric with respect to the origin and it is close with respect to the
sum, under some restrictions ($\m$ and $\n$ have to be linearly
independent, with $|\m|\neq|\n|$ and $\m+\n\in\K_N$).  If one neglects
such restrictions, Lemma \ref{combi} tells us that a set $\enne$ is a
determining set of indices for the cut-off $N$ if it is an algebraic
system of generators for the group $(\z^3,+)$, that is, the smallest
subgroup of $\z^3$ which contains $\enne$ is the whole $\z^3$.

Since by Lemma \ref{combi} it is obviously true that a determining set
of indices, with respect to any cut-off, is a system of generators, one
can ask if the vice-versa is true, that is if each system of generators
is a determining set of indices for a suitable cut-off. We give this
statement in the form of a claim, since in our opinion any proof seems
to be full of technicalities which are not of great interest in this
context.

\begin{claim}
If $\enne$ is an algebraic system of generators for the group $(\z^3,+)$
and $\enne\subset\K_N$, then $\enne$ is a determining set of indices for
the ultraviolet cut-off $N$.
\end{claim}

For the sake of completeness, we give (see Jacobson \cite{Jac}, Theorem
3.8 and Theorem 3.9) a necessary and sufficient condition for a set of
indices $\enne$ to be a system of generators of the whole group $\z^3$.

\begin{theorem}
A set $\enne\subset\z^3$ is a system of generators of $\z^3$ if and only
if the g.c.d.\ of the minors of order $3$ of the matrix $A$ is equal to
$1$, where $A$ is the $k\times 3$ matrix whose rows are the coordinates
of the points of $\enne$ and $k=\#\enne$.
\end{theorem}

The intuitive idea which lets us believe that the claim is true is that
the restrictions given in the statement of property \textsl{(ii)} of
Lemma \ref{combi} can be avoided in the following way.

The restriction about linear independence can be easily avoided by
\textsl{moving aside}: for example if one wants to sum $\m$ with itself,
the best way is to obtain $2\m$ as $\m+\n+\m-\n$, where $\n$ is
linear independent with $\m$.

The restriction about the Euclidean norm (that is, $|\m|\neq|\n|$) can
be avoided, where possible, as in the previous case. Sometimes, as in
the case of the proposition below, this is not possible, since it may
happen that all indices we are allowed to use, have the same Euclidean norm. In such
a case the solution is to \textsl{reach} the index by different paths,
providing with each path a component of the Lie algebra we are dealing
with, in analogy with Lemma \ref{contacci}. This method is probably
peculiar of the dimension three and it does not hold in lower dimensions
(see E and Mattingly \cite{EMa}).

Indeed these tricks are used in the proof of the following proposition,
which states that the \textsl{working example} we talked about in
Section \ref{smain} is a determining set of indices.

\begin{proposition}\label{workex}
Any set $\enne\subset\z^3$ containing the three indices $(1,0,0)$,
$(0,1,0)$ and $(0,0,1)$ is a determining set of indices.
\end{proposition}

\begin{proof}
A careful analysis of the last part of the proof of Lemma \ref{combi}
shows that, if $|\m|=|\n|$, then the Lie brackets $[[F_0,V],W]$,
with $V\in\gU_\m$ and $W\in\gU_\n$, span the two-dimensional subspace
of $\gU_{\m+\n}$ given by
$$
\lambda E\cdot\diff{r_{\m+\n}}+\mu E\cdot\diff{s_{\m+\n}},
$$
where $E$ is the index orthogonal (in $\R^3$) to $\m$ and $\n$.  Hence,
if we sum $(1,0,0)$ and $(0,1,0)$, we obtain the corresponding two
dimensional subspace of $\gU_{(1,1,0)}$. Again a direct computation
shows that such a smaller subspace is indeed sufficient, since if we
combine it with $\gU_{(0,0,1)}$ we obtain the full $\gU_{(1,1,1)}$. Now
we just subtract $(0,0,1)$ from $(1,1,1)$ to obtain the full
$\gU_{(1,1,0)}$ and, in the same way, we can obtain all the indices of
norm $2$. With this set of indices is now easy to obtain, by means of
Lemma \ref{combi} and of the tricks explained above, all the indices in
$\K_N$, whatever is $N$.
\end{proof}

As a consequence of the above proposition is that if $\enne$ is a
determining set of indices for a cut-off $N$, then it is a determining
set of indices for any other cut-off threshold larger than $N$.

\section{The control problem}\label{control}

The section is devoted to the proof of the controllability properties of
the finite dimensional approximations of Navier-Stokes equations. The
first part contains some generalities on polynomial control systems. The
approach and the results are taken from Jurdjevic and Kupka \cite{JuKu}.
In the second part we adapt the proof of a theorem (again of Jurdjevic
and Kupka \cite{JuKu}) to our case. The original theorem applied to
polynomials of odd degree. Polynomials of even degree behave in a
different way, mostly because of the \textsl{obstructions} of the
positive terms. Our case has no obstructions, essentially because of
property \textit{(ii)} of Lemma \ref{contacci}, and the system is
controllable. 

\subsection{Generalities on polynomial control systems}

We consider a system of the form
$$
\dot x=P(x)+\sum_{i=1}^mu_i(t)
$$
where $x\in\R^n$, $b_1$, $b_2,\ldots,b_n$ are fixed vectors in $\R^n$
and $P$ is a polynomial mapping, that is $P=(P_1,\ldots,P_n)$ and each
$P_i$ is a polynomial in the variables $(x_1,\ldots,x_n)$. Let
$Y_1,\ldots,Y_n$ be the constant vector fields assuming respectively value
$b_1,\ldots,b_n$ and let $F$ be the vector fields having the components
of $P$ as its components, and define
$$
\effe=\bigl\{\,F+\sum_{i=1^m}u_iY_i\,\big|\,(u_1,\ldots,u_m)\in\R^m\,\}
$$
We define, for each $x_0\in\R^n$ and $t>0$, the set $A_\effe(x_0,t)$
of states reachable, with a suitable control $u=(u_1,\ldots,u_n)$, from
the initial state $x_0$ in a time smaller than $t$. We define the set
$A^*_\effe(x_0,t)$ of states reachable exactly at time $t$.

Two families of vector fields $\effe_1$, $\effe_2$ are said to be
equivalent if for all $x\in\R^n$ and $t>0$,
$$
\overline{A_{\effe_1}(x,t)}=\overline{A_{\effe_2}(x,t)}.
$$
If $\effe$ is equivalent to $\effe_1$ and to $\effe_2$, then it is
equivalent to $\effe_1\cup\effe_2$. It makes sense then to define the
saturate of $\effe$, denoted by $\sat(\effe)$, which is the union of all
families of vector fields equivalent to $\effe$. Moreover we will call
$\lie(\effe)$ the Lie algebra generated by $\effe$. Finally, the
\textsl{Lie saturate} of $\effe$ is defined as
$\ls(\effe)=\sat(\effe)\cap\lie(\effe)$. In order to obtain
controllability, the Lie saturate should be as large as possible, as
stated by the following theorem.

\begin{theorem}
Let $\effe$ be any family of smooth vector fields and assume that
$\ls(\effe)$ contains $n$ vectors $V_1,\ldots,V_n$ such that the vector
space spanned by them is in $\ls(\effe)$ and for each $x\in\R^n$ the
vectors $V_1(x),\ldots,V_n(x)$ span $\R^n$. Then $A_\effe(x,t)=\R^n$ for
each $x\in\R^n$ and $t>0$.
\end{theorem}

We adapt the conclusions of the theorem to the system that will be
studied in the following section.

\begin{corollary}\label{teorema1}
Let $\effe$ be any family of smooth vector fields and assume that the
constant vector fields of $\ls(\effe)$ span $\R^n$. Then
$A^*_\effe(x,t)=\R^n$ for each $x\in\R^n$ and $t>0$.
\end{corollary}

\begin{proof}
From the previous theorem, $A_\effe(x,t)=\R^n$. Moreover, by Theorem 13,
Chapter 3 of Jurdjevic \cite{Jur} (see also the remarks after Theorem
11, Chapter 5 of \cite{Jur}) it follows that also $A^*_\effe(x,t)=\R^n$.
\end{proof}

In the following, we will need the following two lemmata. The first
lemma permits the \textsl{enlargement} of a family of vector fields by
means of diffeomorphisms. A diffeomorphism $\phi:\R^n\to\R^n$ is a
\textsl{normaliser} of a family $\effe$ if for all $x\in\R^n$ and $t>0$,
$$
\phi\bigl(\overline{A_\effe(\phi^{-1}(x),t)}\bigr)\subset\overline{A_\effe(x,t)},
$$
we will denote by $\norm(\effe)$ the set of all smooth normaliser of
$\effe$.

\begin{lemma}\label{lemma2}
The family $\effe$ is equivalent to
$\bigcup_{\phi\in\norm(\effe)}\{\phi_*(V)\,|\,V\in\effe\}$, where
$\phi_*$ is the differential of $\phi$.
\end{lemma}

The second lemma gives the geometrical structure of the Lie saturate of
a family of vector fields.

\begin{lemma}\label{lemma3}
If $\effe$ is any family of smooth vector fields, then $\effe$ is
equivalent to the closed convex cone generated by $\{\lambda
V\,|\,0\le\lambda\le1,\,V\in\effe\,\}$, where the closure is in the
$C^\infty$ topology on compact sets of $\R^n$.
\end{lemma}

\subsection{Control of the finite dimensional approximations of
Navier-Stokes}

We are able now to prove the controllability property of our
equations. We aim to prove that the control problem
\begin{equation}\label{controleq}
\left\{\begin{array}{l}
\dot r_\ka - F_{r_\ka}(r,s)=q^r_\ka v_\ka^r\\
\dot s_\ka - F_{s_\ka}(r,s)=q^s_\ka v_\ka^s,
\end{array}\right.
\end{equation}
where $F_{r_\ka}$ and $F_{s_\ka}$ are defined in \eqref{fr} and
\eqref{fs}, and the $3\times3$ matrices are defined in \eqref{noiseass},
is controllable, in the sense that for each initial state $(r_I,s_I)\in
U$, for each final state $(r_F,s_F)\in U$ and for each time $T>0$ there
is a family of controls $(v^r_\ka,v^s_\ka)_{\ka\in\enne}$, where $\enne$
is the set of indices corresponding to the non-zero $q_\ka$, such that the
solution corresponding to that control starts at $t=0$ in $(r_I,s_I)$
and arrives in $(r_F,s_F)$ at time $t=T$.

\begin{theorem}\label{controllo}
Assume that the set $\enne$ of non-zero components of the control  is a
determining set of indices, as defined in Section
\ref{determining}. Then system \eqref{controleq} is controllable in the
sense given above.
\end{theorem}

\begin{proof}
First we show that $\gU_\ka\subset\ls(\effe)$ for $\ka\in\enne$. Let
$\lambda\in\R$ and $\ka\in\enne$ and take $Y_\ka\in\gU_\ka$, since
$$
\lambda Y_\ka=\lim_{n\to\infty}\frac1n(F_0+n\lambda Y_\ka),
$$
it follows by Lemma \ref{lemma3} that $Y_\ka\in\ls(\effe)$.

Now we aim to show the following claim: \textsl{if $\m$,
$\n\in\K_N$ are linear independent indices with $|\m|\neq|\n|$ and
$\m+\n\in\K_N$, and if $\gU_\m\oplus\gU_\n\subset\ls(\effe)$, then
also $\gU_{\m+\n}\subset\ls(\effe)$}. If the claim is true, it follows
that each $\gU_\ka$ is contained in $\ls(\effe)$ and, since by the
assumptions $\enne$ is a determining set of indices, Corollary
\ref{teorema1} applies and the proof is ended.

The proof of the claim now follows. From Lemma \ref{combi} we know that
$\gU_{\m+\n}$ is spanned by $[[F_0,V],W]$, where $V\in\gU_\m$ and
$W\in\gU_\n$. From property \textsl{(iii)} of Lemma \ref{contacci}, we
have that
$$
[[F_0,V],W]=\frac12[[F_0,V+W],V+W]
\qquad\text{and}\qquad
-[[F_0,V],W]=\frac12[[F_0,V-W],V-W].
$$
Since by Lemma \ref{lemma3} $\ls(\effe)$ is convex, in order to prove
the claim it is sufficient to show that $\lambda[[F_0,V],V]$ is in
$\ls(\effe)$ for each $\lambda>0$ and $V\in\gU_\m\oplus\gU_\n$.

So, let $\alpha\in\R$ and $V\in\gU_\m\oplus\gU_\n$, then
$\phi(x)=\e^{\alpha V}(x)$ is in $\norm(\effe)$ (see the proof of
Theorem 2 of Jurdjevic and Kupka \cite{JuKu}), so that, by Lemma
\ref{lemma2}, $(\e^{\alpha V})_*(F_0)\in\ls(\effe)$. Now, since the
coefficients of $F_0$ are polynomials of degree $2$ and $V$ is a constant
vector field, it follows that
$$
(\e^{\alpha V})_*(F_0)=I+\alpha[V,F_0]+\frac{\alpha^2}2[V,[V,F_0]],
$$
and so, for each $\lambda>0$,
$$
\lambda[V,[V,F_0]]=\lim_{\alpha\to\infty}\frac\lambda{\alpha^2}(\e^{\alpha
V})_*(F_0)\in\ls(\effe),
$$
since $\ls(\effe)$ is closed. The theorem is proved.
\end{proof}

\section{The exponential convergence}

In this last section we prove Theorem \ref{exporate} as a consequence of
a general result by Meyn and Tweedie \cite{MeTw3} (see Theorem
6.1). Before giving the statement of such theorem, we need to state some
definitions. They will be given in a simplified form, adapted to our
case, while the general statements can be found in the papers by Meyn
and Tweedie \cite{MeTw2}, \cite{MeTw3}.

A nonempty subset $C$ of the state space $U$ is a \textsl{petite} set for a
Markov process with transition probabilities $P_t(\cdot,\cdot)$ if there
are a non-trivial measure $\varphi$ and a probability distribution $a$ on
$(0,\infty)$ such that
$$
\int P_t(x,\cdot)\,a(dt)\ge\varphi\qquad\text{for all }x\in C.
$$

A function $V:U\to\R_+$ is a Lyapunov function for the process if
$V(x)\to\infty$ as $|x|\to\infty$ and there are real constants $c>0$ and
$d$ such that
$$
\gen V(x)\le -cV(x)+d
$$
where $\gen$ is the generator of the diffusion.

The kinetic energy 
$$
V(r,s)=\sum_{k\in\tK}\sum_{i=1}^3({r^i_\ka}^2+{s^i_\ka}^2)
$$
will play the role of the Lyapunov function in our case, as stated by
the following lemma.

\begin{lemma}\label{liapunov}
For each $(r,s)\in U$,
$$
\sum_{\ka\in\tK}\sum_{i=1}^3 (r^i_\ka E_{r^i_\ka}+s^i_\ka E_{s^i_\ka})=0,
$$
where the polynomial $E_{r^i_\ka}$ and $E_{s^i_\ka}$ are respectively
the homogeneous part of degree $2$ of the polynomials $F_{r^i_\ka}$ and
$F_{s^i_\ka}$, defined in \eqref{fr} and \eqref{fs}, and
$$
\gen V(r,s)\le -2\nu V(r,s)+\sigma^2,
$$
where $\gen$ is the generator defined in \eqref{markovgen} and
$\sigma^2$ is the variance of the Brownian motion $B_t$.
\end{lemma}

\begin{proof}
The first property is a far consequence of the famous property of the
non-linear part of Navier-Stokes equations, namely $\int
v\cdot(v\cdot\nabla)v=0$. Indeed
\begin{equation}\label{somme}
\sum_{\ka\in\tK}(r_\ka\cdot E_{r_\ka}+s_\ka\cdot E_{s_\ka})
=\sum_{\ka\in\tK}\sum_{\substack{\h,\el\in\K_N\\
\h+\el=\ka}}\imm[(\ka\cdot u_\h)(u_\el\cdot\overline{u_\ka})]
=\sum_{\ka\in\K_N}u_\ka\cdot E_{u_\ka}
\end{equation}
where $u_\ka=r_\ka+\im s_\ka$ and $E_{u_\ka}$ is the non-linear part in
equation \eqref{finiteq}, namely
$$
E_{u_\ka}=\sum_{\substack{\h,\el\in\K_N\\\h+\el=\ka}}(\ka\cdot
u_\h)\bigl(u_\el-\frac{\ka\cdot u_\el}{|\ka|^2}\ka\bigr).
$$
Finally, the proof that the last sum in \eqref{somme} is equal to $0$ is
just a matter of swapping the two indices $\ka$ and $\el$.

The second property is then an easy consequence of the previous one:
$$
\gen V
=\sum_{\ka\in\tK}\bigl(-2\nu|\ka|^2(r_\ka^2+s_\ka^2)+r_\ka\cdot
E_{r_\ka}+s_\ka\cdot E_{s_\ka}\bigr)+\sigma^2
\le -2\nu V+\sigma^2
$$\end{proof}

Now we are able to prove Theorem \ref{exporate}.

\begin{proof}[Proof of Theorem \ref{exporate}]
From Theorem \ref{main} we know that the Markov process $(r(t),s(t))$ is
strong Feller and irreducible. Using Theorem 3.3 and Theorem 4.1 of Meyn
and Tweedie \cite{MeTw2}, it follows that all compact sets of the state
space $U$ are petite sets. Moreover the previous lemma tells us that the
kinetic energy $V$ is a Lyapunov function. By means of Theorem 6.1 of
Meyn and Tweedie \cite{MeTw3}, we conclude that there are positive
constants $C$ and $\rho$ such that for each $(r_0,s_0)\in U$,
$$
\|P_t((r_0,s_0),\cdot)-\pi\|_f\le C(1+V(r_0,s_0)+\frac{\sigma^2}{2\nu})\e^{-\rho t}.
$$
with $f=1+V$.
\end{proof}

\bibliographystyle{amsplain}

\end{document}